\documentclass{gtart}
\usepackage{amsmath, amssymb, graphicx, epsfig, verbatim, graphs, psfrag}
\usepackage{epic}
\usepackage{amscd}
\usepackage{amssymb}

\newtheorem{thm}{Theorem}[section]

\newtheorem{lem}[thm]{Lemma}
\newtheorem{prop}[thm]{Proposition}
\newtheorem{defn}[thm]{Definition}

\newtheorem{rem}[thm]{Remark}

\numberwithin{equation}{section}

\newcommand{\bfz}{{\mathbb {Z}}}

\newcommand{\bfq}{{\mathbb {Q}}}

\newcommand{\s}{\mathbf s}
\newcommand{\vaa}{\mathbf a}

\renewcommand{\t}{\mathbf t}

\newcommand{\Z}{\mathbb Z}
\newcommand{\N}{\mathbb N}

\newcommand{\R}{\mathbb R}
\newcommand{\bfr}{\mathbb R}

\newcommand{\cpkk}{{\overline {{\mathbb C}{\mathbb P}^2}}}
\newcommand{\cpk}{{\mathbb {CP}}^2}
\newcommand{\cpot}{{\mathbb {CP}}^2\# 5{\overline {{\mathbb C}{\mathbb P}^2}}}
\newcommand{\cphat}{{\mathbb {CP}}^2\# 6{\overline {{\mathbb C}{\mathbb P}^2}}}

\newcommand{\ra}{\rightarrow}

\newcommand{\hf}{{{\widehat {HF}}}}

\newcommand{\frg}{{\mathcal {G}}}

\newcommand{\frw}{{\mathcal {W}}}

\newcommand{\frn}{{\mathcal {N}}}
\newcommand{\farm}{{\mathcal {M}}}

\begin{document}

\title{Symplectic rational blow-down along Seifert fibered 3-manifolds}

\author{David T. Gay}
\address{Department of Mathematics and Applied Mathematics\\ 
University of Cape Town\\
Private Bag X3, Rondebosch 7701\\
South Africa}

\author{Andr\'{a}s I. Stipsicz}
\address{R{\'e}nyi Institute of Mathematics\\
Hungarian Academy of Sciences\\
Re{\'a}ltanoda utca 13--15, \\
H-1053 Budapest, Hungary}

\email{David.Gay@uct.ac.za, stipsicz@math-inst.hu}
\keywords{symplectic rational blow-down, symplectic neighbourhoods, 
rational singularities} 

\begin{abstract}
We verify that the rational blow down schemes along certain
Seifert fibered 3-manifolds found in~\cite{SSW} are, in fact, symplectic
operations.
\end{abstract}
\maketitle

\section{Introduction}

The rational blow-down procedure (introduced by Fintushel and Stern
\cite{FS1} and generalized by Park~\cite{Pratb}) turned out to be one of the
most effective operations in constructing smooth $4$--manifolds with interesting
topological properties, cf.~\cite{FSexo, P, PSS, SS}. Recall that when
performing the rational blow-down operation we simply replace the tubular
neighbourhood of a string of $2$--spheres (intersecting each other according to
the linear plumbing tree with framings specified by the continued fraction
coefficients of $-\frac{p^2}{pq-1}$ for some $p,q$ relatively prime) with a
rational homology disk. The success of this operation might be explained by
the fact that --- as Symington showed~\cite{Sym0, Sym} --- it can be performed
symplectically.  More precisely, if the ambient $4$--manifold is symplectic and
the spheres are symplectic submanifolds intersecting each other orthogonally
then the neighbourhood can be chosen so that the symplectic structure (when
restricted to the complement of this neighbourhood) extends over the glued-in
rational homology disk.  Symington's argument used a beautiful application of
toric geometry by constructing both the right symplectic neighbourhood and the
appropriate symplectic model for the rational homology disk.

It is not hard to list the combinatorial constraints a plumbing tree
must satisfy in order for the proofs of Fintushel--Stern and Park to be
applied. In~\cite{SSW} these contraints were explicitly spelled out
and the family of plumbing trees satisfying these properties has been
determined. The combinatorial conditions are, however, sometimes too weak
to ensure the existence of the rational ho\-mo\-logy disk needed to
perform the geometric operation. In~\cite{SSW} a number of examples
were shown to admit such rational homology disks. (Some of these
examples were already applied by Michalogiorgaki~\cite{michalogiorgaki}
for constructing exotic $4$--manifolds.) Whether these  new examples 
for rational blow-down can be executed in the symplectic category, however,
remained an open question.

In this paper we show that some of the plumbing trees found in~\cite{SSW} can
be, in fact, blown down symplectically.  (For the definition of
\emph{symplectic} rational blow-down see
also~\cite[Definition~1.1]{Sym0}.)  Our main result refers to a class
of plumbing trees which we describe in the following
definition. However, it is worth pointing out here that we are able to
prove the existence of $\omega$--convex neighborhoods (which amounts
to half of our main result) for a much larger class of plumbing trees,
namely all star-shaped negative definite trees.

\begin{defn}\label{d:pelda}
{\rm 
\begin{itemize}
\item 
The plumbing tree given by Figure~\ref{f:wahltype} will be denoted by
$\Gamma _{p,q,r}$ ($p,q,r\geq 0$).
\begin{figure}[ht]
\begin{center}
\setlength{\unitlength}{1mm}
\unitlength=0.7cm
\begin{graph}(10,6)(0,-5)
\graphnodesize{0.2}

 \roundnode{m1}(0,0)
 \roundnode{m2}(1,0)
 \roundnode{m3}(4,0)
 \roundnode{m4}(5,0)
 \roundnode{m5}(6,0)
 \roundnode{m6}(9,0)
 \roundnode{m7}(10,0)
 \roundnode{m8}(5,-1)  
 \roundnode{m9}(5,-4)
 \roundnode{m10}(5,-5)

\edge{m1}{m2}
\edge{m3}{m4}
\edge{m4}{m5}
\edge{m6}{m7}
\edge{m4}{m8}
\edge{m9}{m10}

  \autonodetext{m1}[sw]{{\small $-(p+3)$}}
  \autonodetext{m2}[n]{{\small $-2$}}
  \autonodetext{m3}[n]{{\small $-2$}}
  \autonodetext{m4}[n]{{\small $-4$}}
  \autonodetext{m5}[n]{{\small $-2$}}
  \autonodetext{m6}[n]{{\small $-2$}}
  \autonodetext{m7}[se]{{\small $-(q+3)$}}
  \autonodetext{m8}[w]{{\small $-2$}}
  \autonodetext{m9}[w]{{\small $-2$}}
  \autonodetext{m10}[s]{{\small $-(r+3)$}}
  \autonodetext{m2}[e]{{\Large $\cdots$}}
  \autonodetext{m3}[w]{{\Large $\cdots$}}
  \autonodetext{m5}[e]{{\Large $\cdots$}}
  \autonodetext{m6}[w]{{\Large $\cdots$}}

\freetext(2.5,-0.8)
{$\underset{{\textstyle q}}{\underbrace{\hspace{2.2cm}}}$}

\freetext(7.5,-0.8)
{$\underset{{\textstyle r}}{\underbrace{\hspace{2.2cm}}}$}

\freetext(6,-2.5)
{{\Huge $\rbrace$}}

\freetext(5,-1.3){\Large $.$}
\freetext(5,-1.6){\Large $.$}
\freetext(5,-1.9){\Large $.$}
\freetext(5,-3.1){\Large $.$}
\freetext(5,-3.4){\Large $.$}
\freetext(5,-3.7){\Large $.$}
\freetext(6.5,-2.5){$p$}

\end{graph}
\end{center}
\caption{\quad The graph $\Gamma _{p,q,r}$ in $\frw$ }
\label{f:wahltype}
\end{figure}
We denote the set of these plumbing trees with $\frw$.

\item The plumbing tree of Figure~\ref{f:masik}
will be denoted by $\Delta _{p,q, r}$ ($p\geq 1$ and $q,r\geq 0$).
\begin{figure}[ht]
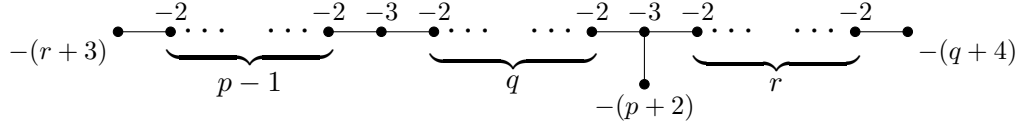

\begin{center}
\setlength{\unitlength}{1mm}
\unitlength=0.7cm
\begin{graph}(6,2)(-1,-1)
\graphnodesize{0.2}

 \roundnode{m12}(-5,0)
 \roundnode{m11}(-4,0)
 \roundnode{m10}(-1,0)
 \roundnode{m1}(0,0)
 \roundnode{m2}(1,0)
 \roundnode{m3}(4,0)
 \roundnode{m4}(5,0)
 \roundnode{m5}(6,0)
 \roundnode{m6}(9,0)
 \roundnode{m7}(10,0)
 \roundnode{m8}(5,-1)  

\edge{m1}{m2}
\edge{m3}{m4}
\edge{m4}{m5}
\edge{m6}{m7}
\edge{m4}{m8}
\edge{m12}{m11}
\edge{m10}{m1}

  \autonodetext{m1}[n]{{\small $-3$}}
  \autonodetext{m2}[n]{{\small $-2$}}
  \autonodetext{m3}[n]{{\small $-2$}}
  \autonodetext{m4}[n]{{\small $-3$}}
  \autonodetext{m5}[n]{{\small $-2$}}
  \autonodetext{m6}[n]{{\small $-2$}}
  \autonodetext{m7}[se]{{\small $-(q+4)$}}
  \autonodetext{m8}[s]{{\small $-(p+2)$}}
 \autonodetext{m10}[n]{{\small $-2$}}
 \autonodetext{m11}[n]{{\small $-2$}}
 \autonodetext{m12}[sw]{{\small $-(r+3)$}} 
 \autonodetext{m2}[e]{{\Large $\cdots$}}
  \autonodetext{m3}[w]{{\Large $\cdots$}}
  \autonodetext{m5}[e]{{\Large $\cdots$}}
  \autonodetext{m6}[w]{{\Large $\cdots$}}
 \autonodetext{m10}[w]{{\Large $\cdots$}}
 \autonodetext{m11}[e]{{\Large $\cdots$}}

\freetext(2.5,-0.8)
{$\underset{{\textstyle q}}{\underbrace{\hspace{2.2cm}}}$}
\freetext(-2.5,-0.8)
{$\underset{{\textstyle p-1}}{\underbrace{\hspace{2.2cm}}}$}

\freetext(7.5,-0.8)
{$\underset{{\textstyle r}}{\underbrace{\hspace{2.2cm}}}$}
\end{graph}
\end{center}
\caption{\quad The graph $\Delta _{p,q,r}$ for $p\geq 1$ and $q,r\geq 0$}
\label{f:masik}
\end{figure}
The slight modification of the graph $\Delta _{p,q,r}$ when $p=0$ is shown
in Figure~\ref{f:specmasik}. The set of graphs $\Delta _{p,q,r}$ with
$p,q,r\geq 0$ will be denoted by $\frn$.

\begin{figure}[ht]
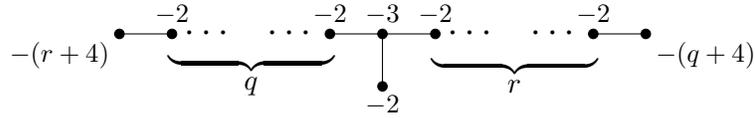

\begin{center}
\setlength{\unitlength}{1mm}
\unitlength=0.7cm
\begin{graph}(10,4)(0,-1)
\graphnodesize{0.2}

 \roundnode{m1}(0,0)
 \roundnode{m2}(1,0)
 \roundnode{m3}(4,0)
 \roundnode{m4}(5,0)
 \roundnode{m5}(6,0)
 \roundnode{m6}(9,0)
 \roundnode{m7}(10,0)
 \roundnode{m8}(5,-1)  

\edge{m1}{m2}
\edge{m3}{m4}
\edge{m4}{m5}
\edge{m6}{m7}
\edge{m4}{m8}

  \autonodetext{m1}[sw]{{\small $-(r+4)$}}
  \autonodetext{m2}[n]{{\small $-2$}}
  \autonodetext{m3}[n]{{\small $-2$}}
  \autonodetext{m4}[n]{{\small $-3$}}
  \autonodetext{m5}[n]{{\small $-2$}}
  \autonodetext{m6}[n]{{\small $-2$}}
  \autonodetext{m7}[se]{{\small $-(q+4)$}}
  \autonodetext{m8}[s]{{\small $-2$}}
  \autonodetext{m2}[e]{{\Large $\cdots$}}
  \autonodetext{m3}[w]{{\Large $\cdots$}}
  \autonodetext{m5}[e]{{\Large $\cdots$}}
  \autonodetext{m6}[w]{{\Large $\cdots$}}

\freetext(2.5,-0.8)
{$\underset{{\textstyle q}}{\underbrace{\hspace{2.2cm}}}$}

\freetext(7.5,-0.8)
{$\underset{{\textstyle r}}{\underbrace{\hspace{2.2cm}}}$}
\end{graph}
\end{center}
\caption{\quad The graph $\Delta _{0,q,r}$}
\label{f:specmasik}
\end{figure}

\item The plumbing graph of Figure~\ref{f:lambda} is $\Lambda _{p,q,r}$
with $p, r\geq 1$ and $q\geq 0$.
\begin{figure}[ht]
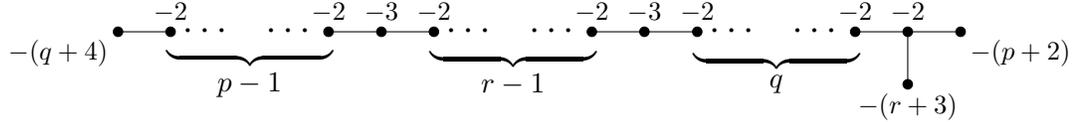

\begin{center}
\setlength{\unitlength}{1mm}
\unitlength=0.7cm
\begin{graph}(-6,2)(0,-1)
\graphnodesize{0.2}

 \roundnode{m15}(-10,0)
 \roundnode{m14}(-9,0)
 \roundnode{m13}(-6,0)
 \roundnode{m12}(-5,0)
 \roundnode{m11}(-4,0)
 \roundnode{m10}(-1,0)
 \roundnode{m1}(0,0)
 \roundnode{m2}(1,0)
 \roundnode{m3}(4,0)
 \roundnode{m4}(5,0)
 \roundnode{m5}(6,0)
  \roundnode{m8}(5,-1)  

\edge{m1}{m2}
\edge{m3}{m4}
\edge{m4}{m5}
\edge{m4}{m8}
\edge{m12}{m11}
\edge{m10}{m1}
\edge{m15}{m14}
\edge{m13}{m12}

  \autonodetext{m1}[n]{{\small $-3$}}
  \autonodetext{m2}[n]{{\small $-2$}}
  \autonodetext{m3}[n]{{\small $-2$}}
  \autonodetext{m4}[n]{{\small $-2$}}
  \autonodetext{m5}[se]{{\small $-(p+2)$}}
  \autonodetext{m8}[s]{{\small $-(r+3)$}}
 \autonodetext{m10}[n]{{\small $-2$}}
 \autonodetext{m11}[n]{{\small $-2$}}
 \autonodetext{m12}[n]{{\small $-3$}} 
  \autonodetext{m13}[n]{{\small $-2$}}
 \autonodetext{m14}[n]{{\small $-2$}}  
 \autonodetext{m15}[sw]{{\small $-(q+4)$}} 
\autonodetext{m2}[e]{{\Large $\cdots$}}
\autonodetext{m3}[w]{{\Large $\cdots$}}
  \autonodetext{m10}[w]{{\Large $\cdots$}}
 \autonodetext{m11}[e]{{\Large $\cdots$}}
\autonodetext{m14}[e]{{\Large $\cdots$}}
\autonodetext{m13}[w]{{\Large $\cdots$}}

\freetext(2.5,-0.8)
{$\underset{{\textstyle q}}{\underbrace{\hspace{2.2cm}}}$}
\freetext(-2.5,-0.8)
{$\underset{{\textstyle r-1}}{\underbrace{\hspace{2.2cm}}}$}

\freetext(-7.5,-0.8)
{$\underset{{\textstyle p-1}}{\underbrace{\hspace{2.2cm}}}$}
\end{graph}
\end{center}
\caption{\quad The graph $\Lambda _{p,q,r}$ for $p,r\geq 1$ and $q\geq 0$}
\label{f:lambda}
\end{figure}
Modifications of these graphs for $p=0, r\geq 1$ and $p\geq 1, r=0$ and
finally for $p=r=0$ are shown by Figures~\ref{f:lambdaspec1},
\ref{f:lambdaspec2} and by Figure~\ref{f:cs}.  The set of graphs $\Lambda
_{p,q,r}$ with $ p,q,r\geq 0$ will be denoted by $\farm$.

\begin{figure}[ht]
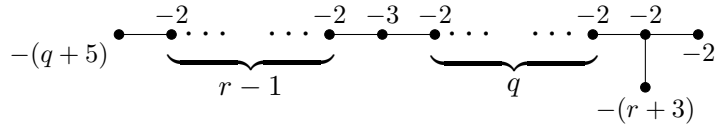

\begin{center}
\setlength{\unitlength}{1mm}
\unitlength=0.7cm
\begin{graph}(-2,2)(0,-1)
\graphnodesize{0.2}

 \roundnode{m12}(-5,0)
 \roundnode{m11}(-4,0)
 \roundnode{m10}(-1,0)
 \roundnode{m1}(0,0)
 \roundnode{m2}(1,0)
 \roundnode{m3}(4,0)
 \roundnode{m4}(5,0)
 \roundnode{m5}(6,0)
  \roundnode{m8}(5,-1)  

\edge{m1}{m2}
\edge{m3}{m4}
\edge{m4}{m5}
\edge{m4}{m8}
\edge{m12}{m11}
\edge{m10}{m1}

  \autonodetext{m1}[n]{{\small $-3$}}
  \autonodetext{m2}[n]{{\small $-2$}}
  \autonodetext{m3}[n]{{\small $-2$}}
  \autonodetext{m4}[n]{{\small $-2$}}
  \autonodetext{m5}[s]{{\small $-2$}}
  \autonodetext{m8}[s]{{\small $-(r+3)$}}
 \autonodetext{m10}[n]{{\small $-2$}}
 \autonodetext{m11}[n]{{\small $-2$}}
 \autonodetext{m12}[sw]{{\small $-(q+5)$}} 
\autonodetext{m2}[e]{{\Large $\cdots$}}
\autonodetext{m3}[w]{{\Large $\cdots$}}
 \autonodetext{m10}[w]{{\Large $\cdots$}}
 \autonodetext{m11}[e]{{\Large $\cdots$}}

\freetext(2.5,-0.8)
{$\underset{{\textstyle q}}{\underbrace{\hspace{2.2cm}}}$}
\freetext(-2.5,-0.8)
{$\underset{{\textstyle r-1}}{\underbrace{\hspace{2.2cm}}}$}

\end{graph}
\end{center}
\caption{\quad The graph $\Lambda _{0,q,r}$ for $r\geq 1$ and $q\geq 0$}
\label{f:lambdaspec1}
\end{figure}

\begin{figure}[ht]
\begin{center}
\setlength{\unitlength}{1mm}
\unitlength=0.7cm
\begin{graph}(-2,2)(0,-1)
\graphnodesize{0.2}

 \roundnode{m12}(-5,0)
 \roundnode{m11}(-4,0)
 \roundnode{m10}(-1,0)
 \roundnode{m1}(0,0)
 \roundnode{m2}(1,0)
 \roundnode{m3}(4,0)
 \roundnode{m4}(5,0)
 \roundnode{m5}(6,0)
  \roundnode{m8}(5,-1)  

\edge{m1}{m2}
\edge{m3}{m4}
\edge{m4}{m5}
\edge{m4}{m8}
\edge{m12}{m11}
\edge{m10}{m1}

  \autonodetext{m1}[n]{{\small $-4$}}
  \autonodetext{m2}[n]{{\small $-2$}}
  \autonodetext{m3}[n]{{\small $-2$}}
  \autonodetext{m4}[n]{{\small $-2$}}
  \autonodetext{m5}[se]{{\small $-(p+2)$}}
  \autonodetext{m8}[s]{{\small $-3$}}
 \autonodetext{m10}[n]{{\small $-2$}}
 \autonodetext{m11}[n]{{\small $-2$}}
 \autonodetext{m12}[sw]{{\small $-(q+4)$}} 
\autonodetext{m2}[e]{{\Large $\cdots$}}
\autonodetext{m3}[w]{{\Large $\cdots$}}
 \autonodetext{m10}[w]{{\Large $\cdots$}}
 \autonodetext{m11}[e]{{\Large $\cdots$}}

\freetext(2.5,-0.8)
{$\underset{{\textstyle q}}{\underbrace{\hspace{2.2cm}}}$}
\freetext(-2.5,-0.8)
{$\underset{{\textstyle p-1}}{\underbrace{\hspace{2.2cm}}}$}

\end{graph}
\end{center}
\caption{\quad The graph $\Lambda _{p,q,0}$ for $p\geq 1$ and $q\geq 0$}
\label{f:lambdaspec2}
\end{figure}

\begin{figure}[ht]
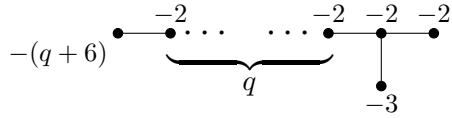

\begin{center}
\setlength{\unitlength}{1mm}
\unitlength=0.7cm
\begin{graph}(1,2)(2,-1)
\graphnodesize{0.2}

 \roundnode{m1}(0,0)
 \roundnode{m2}(1,0)
 \roundnode{m3}(4,0)
 \roundnode{m4}(5,0)
 \roundnode{m5}(6,0)
  \roundnode{m8}(5,-1)  

\edge{m1}{m2}
\edge{m3}{m4}
\edge{m4}{m5}
\edge{m4}{m8}

  \autonodetext{m1}[sw]{{\small $-(q+6)$}}
  \autonodetext{m2}[n]{{\small $-2$}}
  \autonodetext{m3}[n]{{\small $-2$}}
  \autonodetext{m4}[n]{{\small $-2$}}
  \autonodetext{m5}[n]{{\small $-2$}}
  \autonodetext{m8}[s]{{\small $-3$}}
  \autonodetext{m2}[e]{{\Large $\cdots$}}
  \autonodetext{m3}[w]{{\Large $\cdots$}}

\freetext(2.5,-0.8)
{$\underset{{\textstyle q}}{\underbrace{\hspace{2.2cm}}}$}

\end{graph}
\end{center}
\caption{\quad The graph $\Lambda _{0,q,0}$ for $q\geq 0$}
\label{f:cs}
\end{figure}
\item Finally let us denote the union $\frw\cup \frn \cup \farm$ simply
by $\frg$.
\end{itemize}
}
\end{defn}

\begin{thm}\label{t:main}
  Suppose that $(X, \omega )$ is a given symplectic $4$--manifold and the
  symplectic spheres $S_i\subset (X, \omega )$ intersect each other
  $\omega$--orthogonally and according to the plumbing tree $\Gamma$, which is
  an element of $\frg$.  Then there is a rational homology disk $B_{\Gamma}$
  and a tubular neighbourhood $S_{\Gamma }$ of $\cup S_i\subset X$ such that
  the symplectic form $\omega $ extends from $X-{\mbox {int}}S_{\Gamma }$ to
  $X_{\Gamma }= (X-{\mbox {int}}S_{\Gamma })\cup B_{\Gamma}$. In short, graphs
  in $\frg$ can be symplectically blown down.
\end{thm}

The proof relies on the gluing theorem of symplectic manifolds along
$\omega$--convex (and $\omega$--concave) boundaries, as explained
in~\cite{etny}. (Recall that ``$\omega$--convex'', resp.
``$\omega$--concave'', is synonomous with ``strongly convex'', resp.
``strongly concave''.) According to this scheme, a strong convex symplectic
filling of the contact $3$--manifold $(M_1, \xi _1)$ can be symplectically
glued to a strong concave symplectic filling of $(M_2, \xi _2)$ provided there
is a contactomorphism between $(M_1, \xi _1)$ and $(M_2, \xi _2)$. Therefore
one way of proving Theorem~\ref{t:main} is to show that
\begin{enumerate}
\item the symplectic spheres admit a neighbourhood such that the complement 
of it is a strong concave filling of the boundary contact $3$--manifold,
\item the rational homology disk $B_{\Gamma}$ admits a symplectic structure which is a
strong convex filling of its boundary (equipped with some contact structure),
and finally
\item the desired gluing contactomorphism exists.
\end{enumerate}

There are many possibilities for verifying (1) above. Since the
plumbing graph under examination is negative definite, a classical
theorem of Grauert \cite{Grau} shows the existence of a $J$--convex
neighbourhood of the spheres, which in this di\-men\-sion shows that
the complement is a weak concave filling of its boundary. The concave
analogue of Eliashberg's deformation argument~\cite{Eli} (stating that
a weak convex filling of a rational homology $3$--sphere can be
perturbed to be a strong filling, cf.~\cite{OO}) would complete the
argument --- but no concave analogue of the ``convex'' theorem of
Eliashberg mentioned above is available at the moment.  We rather use
an explicit way of constructing symplectic structures on the plumbing
$4$--manifold $M_{\Gamma}$ (associated to the plumbing graph $\Gamma$)
and invoke standard neighbourhood theorems to find the right
$\omega$--convex neighbourhoods.

For (2) we only need to note that the rational homology disks for the
plumbing graphs of $\frg$ are given as deformations of appropriate
surface singularities (see~\cite{SSW}), therefore are equipped with
Stein structures. (Alternatively, $B_{\Gamma}$ can be given by
Weinstein handle attachments, according to the Kirby diagram described
in \cite{SSW}.)

Finally, the existence of the desired contactomorphism (listed under
(3) above) will be shown by using the classification of tight contact
structures on the boundary $3$--manifolds.

{\bf {Acknowledgements:}} This material is based upon work supported by the
National Research Foundation of South Africa under Grant number 62124
(Hungarian Grant Number ZA-15/2006). The second author was also partially
supported by OTKA T49449 and by the EU Marie Curie TOK grant BudAlgGeo. We
would also like to thank the referee for helpful comments which have improved
the exposition.

\section{Convex neighbourhoods of certain plumbings}

We prove a slightly stronger result than is needed in the proof of our
main result. More precisely, we will consider a more general class of
graphs than $\frg$.

\begin{thm}\label{t:nbhood}
Suppose that $(X, \omega )$ is a given symplectic $4$--manifold, and the
symplectic spheres $S_i$ intersect each other $\omega$--orthogonally
and according to a star-shaped graph $\Gamma$ with negative definite
intersection form. Then $\cup S_i$ admits an $\omega$--convex
neighbourhood $S_{\Gamma}$.
\end{thm}

The proof of Theorem~\ref{t:nbhood} will easily follow from the main result of
this section below. (As customary, we will denote the positive half line $(0,
\infty )$ by $\bfr ^+$.)
\begin{thm}\label{t:toric}
Suppose that $\Gamma$ is a star-shaped negative definite graph and
$a_i\in \bfr ^+$ ($i=1, \ldots , \vert \Gamma \vert $) are given. Then
there is a symplectic structure $\omega _{\vaa}$ on the plumbing
$4$--manifold $M_{\Gamma }$ associated to the graph $\Gamma \in \frg$
such that
\begin{itemize}
\item the spheres $P_i$ corresponding to the vertices of the plumbing
graph are symplectic and $\omega$--orthogonal,
\item $\int _{P_i}\omega _{\vaa}=a_i$ and  
\item any neighbourhood of $\cup P_i$ contains an $\omega _{\vaa}$--convex
neighbourhood.
\end{itemize}
\end{thm}

\begin{proof}[Proof of Theorem~\ref{t:nbhood} from Theorem~\ref{t:toric}]
  Suppose that $(X, \omega )$ with $S_i\subset X$ intersecting each
  other according to $\Gamma$ is given and take $a_i =\int _{S_i}
  \omega$.  Consider the plumbing $4$--manifold $M_{\Gamma}$
  corresponding to $\Gamma$ and equip $M_{\Gamma}$ with the symplectic
  structure $\omega _{\vaa}$ provided by Theorem~\ref{t:toric}. By
  Moser's Theorem $\cup S_i \subset (X, \omega )$ and $\cup P_i
  \subset (M_{\Gamma }, \omega _{\vaa})$ admit symplectomorphic
  neighbourhoods (cf. also~\cite[Proposition~3.5]{Sym0}), which, by
  the third conclusion of Theorem~\ref{t:toric} concludes the proof.
\end{proof}

The idea of the proof of Theorem~\ref{t:toric} is to construct the
desired neighbourhood for each of the legs independently, using toric
techniques as in~\cite{Sym0}, then to construct the desired
neighbourhood for the central vertex, and then to glue the pieces
together carefully. These techniques will certainly generalize to
larger classes of graphs, but here we focus on star-shaped
graphs. Also, there should be a way to state these results in terms of
embedded graphs in $\bfr ^3$ with edges of rational slope, inspired by
the ``tropical'' ideas of Mikhalkin, but in these simple cases it is
easier to avoid this perspective.

It is convenient to  work with $5$--tuples of the form
$(X,\omega,C,f,V)$ where:
\begin{enumerate}
\item $X$ is a $4$-manifold, 
\item $\omega$ is a symplectic form on $X$,
\item $C$ is a collection of symplectically embedded
surfaces in $X$ intersecting each other $\omega$--orthogonally, 
\item $f$ is a smooth
function $f\colon X \ra [0,\infty)$ such that $f^{-1}(0) = C$ and $f$ has no
critical values in $(0,\infty)$, and 
\item $V$ is a Liouville vector field defined on $X \setminus C$ such
  that $df(V) > 0$.
\end{enumerate}
If we can produce such a $5$--tuple with $C$ being a collection of
spheres intersecting according to the given graph $\Gamma$ with
symplectic areas $a_i$, then we will have proved
Theorem~\ref{t:toric}. In our proof we will construct one such
$5$--tuple for each leg of $\Gamma$ and one for the central vertex, and
then glue them together.

Let us call such a $5$--tuple a ``neighbourhood $5$--tuple''.
Theorem~\ref{t:toric} follows from a much more technical statement:

\begin{prop} \label{p:technicaltoric}
Suppose that, for some $m \in \N$ and some $n_1, \ldots, n_m \in \N$,
we are given points in the upper half plane $\{P_{i,j} =
(x_{i,j},y_{i,j}) \in \R^2, y_{i,j} > 0 \mid 1 \leq i \leq m, 0 \leq j
\leq n_i+2\}$, with $P_{i,j} \neq P_{i,j+1}$. For each $i,j$ with $1
\leq i \leq m$ and $0 \leq j \leq n_i+1$, let $L_{i,j}$ be the line
containing $P_{i,j}$ and $P_{i,j+1}$ and let $E_{i,j}$ be the compact
oriented line segment from $P_{i,j}$ to $P_{i,j+1}$. Now require that
each $L_{i,j}$ has rational or infinite slope, so that there is a
well-defined primitive integral tangent vector $\tau_{i,j} =
(u_{i,j},v_{i,j})^T \in \Z^2$ parallel to $E_{i,j}$ pointing from
$P_{i,j}$ to $P_{i,j+1}$. (Here ``primitive'' means that
$\gcd(u_{i,j},v_{i,j})=1$ with $\gcd(1,0)$ defined to be $1$.) Record
the positive real numbers $\lambda_{i,j}$ such that
the vector from $P_{i,j}$ to $P_{i,j+1}$ is equal to $\lambda_{i,j}
\tau_{i,j}$. Suppose there is a fixed constant $y_0 > 0$ such that
$y_{i,0} = y_0$ for all $i = 1, \ldots, m$ and suppose that $x_{1,0} +
x_{2,0} + \ldots + x_{m,0} > 0$. Suppose furthermore that the
following conditions are satisfied for each $i = 1, \ldots, m$:
\begin{enumerate}
\item $\tau_{i,0} = (1,0)^T$; i.e. $x_{i,0} < x_{i,1}$ and $y_{i,0} =
  y_{i,1}$.
%\item There is a fixed constant $y_0 > 0$ such that $y_{i,0} = y_0$.
\item $L_{i,n_i+1}$ passes through the origin $O = (0,0)$;
  i.e. $(x_{i,n_i+2},y_{i,n_i+2}) = \rho (x_{i,n_i+1},y_{i,n_i+1})$
  for some $\rho \in \R^+$.
\item For each $j=0, \ldots, n_i$, $\det(\tau_{i,j},\tau_{i,j+1}) =
+1$.
\item For each $j=0, \ldots, n_i$, $\det(\tau_{i,j}, P_{i,j}^T) > 0$.
\end{enumerate}
Then there exists a neighbourhood $5$--tuple $(X,\omega,C,f,V)$ such
that $C$ is a collection of symplectically embedded $2$--spheres $S_0,
S_{i,j}, 1\leq i \leq m, 1 \leq j \leq n_i$ intersecting
$\omega$--orthogonally according to an $m$--legged star-shaped graph with
$S_0$ the central vertex intersecting each $S_{i,1}$ and with
$S_{i,1}, \ldots, S_{i,n_i}$ the $i$'th leg. Furthermore, the areas
and self-intersections are given as follows:
\begin{enumerate}
\item For each $i = 1, \ldots, m$ and each $j = 1, \ldots, n_i$,
  $S_{i,j} \cdot S_{i,j} = - \det(\tau_{i,j-1},\tau_{i,j+1})$, and the
  area of $S_{i,j}$ is $2 \pi \lambda_{i,j}$.
\item Noting that, for each $i=1, \ldots, m$, $\tau_{i,1} =
  (u_{i,1},1)^T$, we have that $S_0 \cdot S_0 = u_{1,1} + u_{2,1} + \ldots +
  u_{m,1}$.
\item The area of $S_0$ is $2\pi(x_{1,1} + x_{2,1} + \ldots +
x_{m,1})$.
\end{enumerate}
\end{prop}

\begin{proof}

  Throughout this proof, identify $S^1$ with $\R / 2\pi \Z$ through the $\R /
  2\pi \Z$--valued coordinate function $q$ (possibly decorated with
  subscripts); i.e. $\int_{S^1} dq = 2\pi$.

To build the neighbourhood $5$--tuples for the legs, we give a detailed
description of the toric setup as we need it and the basic results
that come from the techniques in~\cite{Sym0, Sym}. (We are essentially
ignoring the global group action point of view that is fundamental to
traditional toric geometry, and taking the more local topological
point of view developed by Symington. Also, here we make essential use
of the correspondence between radial vector fields on images of
moment maps and Liouville vector fields on the domains of moment maps;
this is not standard in the toric geometry literature but is critical
to Symington's work and ours.)

For each $i=1, \ldots, m$, the following is the construction of a
neighbourhood $5$--tuple $(X_i,\omega_i,C_i,f_i,V_i)$ for the $i$'th
leg. Here $C_i$ will be $D_{i,0} \cup S_{i,1} \cup S_{i,2} \cup \ldots
\cup S_{i,n}$ where $D_{i,0}$ is a small open disk neighbourhood in
$S_0$ of $S_0 \cap S_{i,1}$. An example is illustrated in
Figure~\ref{f:oneleg}.
\begin{figure}[ht]
\begin{center}
\psfrag{Ri}{$R_i$}
\psfrag{Pi0}{$P_{i,0}$}
\psfrag{Pi1}{$P_{i,1}$}
\psfrag{Pi2}{$P_{i,2}$}
\psfrag{Pi3}{$P_{i,3}$}
\psfrag{Pi4}{$P_{i,4}$}
\psfrag{Ei0}{$E_{i,0}$}
\psfrag{Ei1}{$E_{i,1}$}
\psfrag{Ei2}{$E_{i,2}$}
\psfrag{Ei3}{$E_{i,3}$}
\psfrag{f-1e}{$g_i^{-1}(\epsilon)$}
\includegraphics[width=12cm]{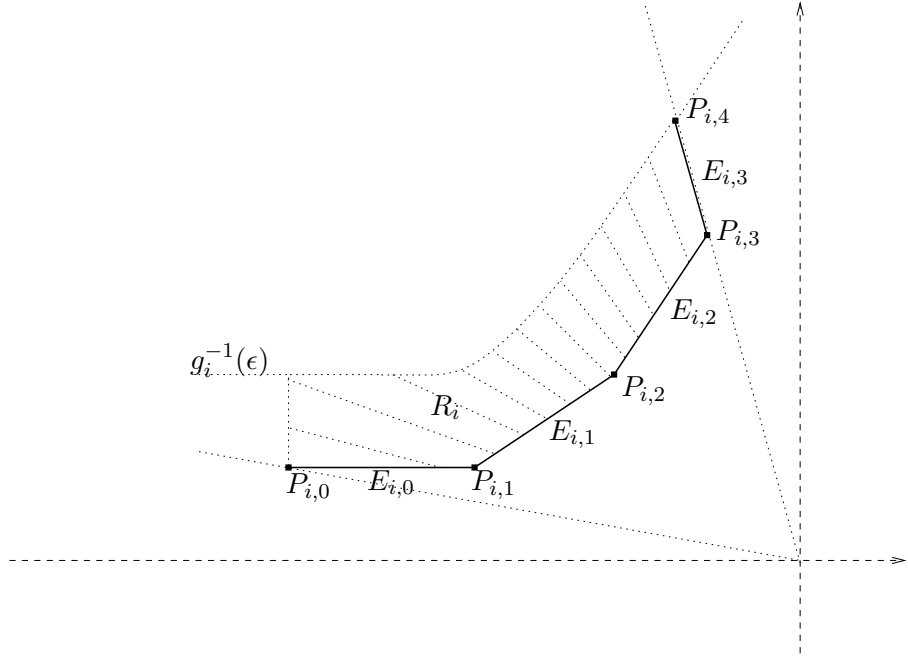}
\end{center}
\caption{An example of a two-dimensional ``moment-map image'' of the
neighbourhood for one leg of a star-shaped plumbing graph. Here $n_i$ =
2.}
\label{f:oneleg}
\end{figure}

For each $j=0,\ldots, n_i+1$, let $H_{i,j}$ be the closed half-plane
with $\partial H_{i,j} = L_{i,j}$ and not containing the origin $O$,
and consider the convex region $H_i=H_{i,0} \cap \ldots \cap
H_{i,n_i+1}$. Choose a smooth function $g_i\colon  H_i \ra [0,\infty)$ with
no critical values in $(0,\infty)$ such that:
\begin{enumerate}
\item $g_i^{-1}(0) = (L_{i,0} \cup L_{i,1} \cup \ldots \cup L_{i,n_i})
\cap H_i$.
\item For some small $\epsilon > 0$ (smaller than $x_{i,1}-x_{i,0}$),
  $g_i$ restricted to $(-\infty,x_{i,0}+\epsilon] \times
  [y_{i,0},y_{i,0}+\epsilon]$ is given by $g_i(x,y) = y-y_{i,0}$.
\item On a small neighbourhood of $L_{i,n_i+1} \cap H_i$, $g_i(x,y)$ is
  equal to the Euclidean distance from $(x,y)$ to $L_{i,n_i}$.
\item The level sets of $g_i$ are transverse to the radial vector field
  radiating out from the origin; i.e. $dg_i(x \partial_x + y \partial_y)
  > 0$ whenever $dg_i \neq 0$.
\end{enumerate}
Finally consider the region $R_i = g_i^{-1}[0,\epsilon) \setminus
((-\infty,x_{i,0}] \times \R)$, and relabel things so that $E_{i,0}
\cap R_i$ is now called $E_{i,0}$ and $E_{i,n_i+1} \cap R_i$ is now
called $E_{i,n_i+1}$. (These are the initial and final edges of the
polygonal boundary of $R_i$, and they are bounded but noncompact,
while the other edges $E_{i,1}, \ldots, E_{i,n}$ are compact.)  Also
let $P_{i,n_i+2}$ now denote the new terminal point of $E_{i,n_i+1}$,
so that $E_{i,n_i+1}$ is the half-open interval from $P_{i,n_i+1}$ to
$P_{i,n_i+2}$, containing $P_{i,n_i+1}$ and not $P_{i,n_i+2}$. This
$2$--dimensional data is the ``template'' for our neighbourhood
$5$--tuple, which we now describe in some detail.

Let $Y_i = \{(p_1,q_1,p_2,q_2) \in \R \times S^1 \times \R \times S^1 \mid
(p_1,p_2) \in R_i \}$ with symplectic form $\eta_i = dp_1 \wedge dq_1 + dp_2
\wedge dq_2$ and let $\pi_i\colon Y_i \ra R_i$ be the projection
$(p_1,q_1,p_2,q_2) \mapsto (p_1,p_2)$. For each edge $E_{i,j}$ let $Y_{i,j} =
[0,\lambda_{i,j}] \times S^1$ with coordinates $(p,q)$ and symplectic form
$\eta_{i,j} = dp \wedge dq$ and let $\pi_{i,j}\colon Y_{i,j} \ra E_{i,j}$ be
the projection defined by $(p,q) \mapsto (x_{i,j} + p u_{i,j}, y_{i,j} + p
v_{i,j})$. Also, for each $j$, consider the map $\psi_{i,j}\colon
\pi_i^{-1}(E_{i,j}) \ra Y_{i,j}$ defined by $(p_1,q_1,p_2,q_2) \mapsto (p,q)$
where $p$ is such that $(p_1,p_2) = (x_{i,j} + p u_{i,j}, y_{i,j} + p
v_{i,j})$ and $q = u_{i,j} q_1 + v_{i,j} q_2$. In other words, the $S^1 \times
S^1$ in $Y_i$ above a point $(p_1,p_2) \in E_{i,j}$ collapses to the $S^1$ in
$Y_{i,j}$ above $(p_1,p_2)$, with the $(-v_{i,j},u_{i,j})$--curves in $S^1
\times S^1$ being the curves that collapse to points. Note that on
$\pi_i^{-1}(E_{i,j})$ we have $\pi_i = \pi_{i,j} \circ \psi_{i,j}$.
Then there exists a unique symplectic $4$--manifold $(X_i,\omega_i)$
equipped with the following smooth maps:
\begin{enumerate}
\item A map $\mu_i\colon X_i \ra \R^2$ with $\mu_i(X_i) = R_i$. (This is the
  ``moment map'' for a Hamiltonian torus action on $(X_i,\omega_i)$,
  after identifying $\R^2$ with the dual of the Lie algebra of $S^1
  \times S^1$.)
\item A map $\phi_i\colon Y_i \ra X_i$ such that $\pi_i = \mu_i \circ
  \phi_i$ and such that $\phi_i$ is a symplectomorphism from
  $(\pi_i^{-1}(R_i \setminus \partial R_i),\eta_i)$ to
  $(\mu_i^{-1}(R_i \setminus \partial R_i),\omega_i)$. We think of
  this as giving us coordinates $(p_1,q_1,p_2,q_2)$ on $X_i$ with
  respect to which $\omega_i = dp_1 \wedge dq_1 + dp_2 \wedge dq_2$,
  except that the $S^1$--valued coordinates $(q_1,q_2)$ degenerate
  along $\mu_i^{-1}(\partial R_i)$ in various ways determined by the
  next set of maps.
\item A map $\phi_{i,j}\colon Y_{i,j} \ra \mu_i^{-1}(E_{i,j})$, for each $j
  = 0, \ldots, n_i+1$, such that $\pi_{i,j} = \mu_i \circ \phi_{i,j}$
  and such that $\phi_{i,j}$ restricts to a symplectomorphism from
  $(\pi_{i,j}^{-1}(E_{i,j} \setminus \partial E_{i,j}), \eta_{i,j})$
  to $(\mu_i^{-1}(E_{i,j} \setminus \partial E_{i,j}), \omega_i)$.
\end{enumerate}
In addition, for each point $P_{i,j} \in R_i$, $\mu_i^{-1}(P_{i,j})$
is a single point. 
The following are consequences of this fact and the relationships
amongst these maps:
\begin{enumerate}
\item For each point $P \in R_i \setminus \partial R_i$, $\mu_i^{-1}(P)$ is
  a torus.
\item For each point $P \in E_{i,j} \setminus \partial E_{i,j}$,
  $\mu_i^{-1}(P)$ is a circle.
\item For each compact edge $E_{i,j}$, $j=1, \ldots, n_i$,
  $\mu_i^{-1}(E_{i,j})$ is a symplectically embedded $2$--sphere
  $S_{i,j}$ with symplectic area $2\pi \lambda_{i,j}$.
\item For each of the two noncompact edges $E_{i,0}$ and $E_{i,n_i+1}$,
  $D_{i,j} = \mu_i^{-1}(E_{i,j})$ is a symplectically embedded open disk with
  symplectic area $2\pi \lambda_{i,j}$ ($j=0$ or $n_i+1$).
\item Each point $\mu_i^{-1}(P_{i,j})$, for $j=1, \ldots, n_i+1$, is a point
  of $\omega$--orthogonal intersection between $\mu_i^{-1}(E_{i,j-1})$
  and $\mu_i^{-1}(E_{i,j})$.
\item Given an embedded arc $\gamma\colon [0,1] \ra R_i$ meeting $\partial
  R_i$ transversely at one point $\gamma(1)$ in the interior of an
  edge $E_{i,j}$, the submanifold $\mu_i^{-1}(\gamma([0,1]))$ is a
  solid torus naturally parameterized using the coordinates
  $(t,q_1,q_2)$ as $[0,1] \times S^1 \times S^1$, with $\{1\} \times S^1
  \times S^1$ collapsed to $S^1$ so that the $(-v_{i,j},u_{i,j})$
  curves collapse to points.
\item Each $2$--sphere $S_{i,j} = \mu_i^{-1}(E_{i,j})$, for $j=1,
  \ldots, n_i$, has self-intersection $S_{i,j} \cdot S_{i,j} =
  -\det(\tau_{i,j-1}, \tau_{i,j+1})$. This can be seen by using the
  preceding point to understand the topology of the $S^1$--bundle over
  $S^2$ which is the boundary of a tubular neighborhood of
  $S_{i,j}$. This bundle is the inverse image under $\mu_i$ of a
  parallel copy of $E_{i,j}$, translated into the interior of $R_i$
  and extended until it intersects the two adjacent edges $E_{i,j-1}$
  and $E_{i,j+1}$. Thus the bundle is the result of Dehn filling the
  two boundary components of $[0,1] \times S^1 \times S^1$ according
  to the rule in the preceding point, which gives an Euler class of
  $-\det(\tau_{i,j-1}, \tau_{i,j+1})$.
\end{enumerate}
Furthermore, focusing attention on the end
$\mu_i^{-1}((x_{i,0},x_{i,0}+\epsilon) \times
[y_{i,0},y_{i,0}+\epsilon)) $, which is diffeomorphic to
$(x_{i,0},x_{i,0}+\epsilon) \times S^1 \times D$, where $D$ is an open
disk in $\R^2$ of an appropriate radius, everything can be written
explicitly as follows, using coordinates $t \in
(x_{i,0},x_{i,0}+\epsilon)$, $\alpha \in S^1$ and polar coordinates
$(r,\theta) \in D$:
\begin{enumerate}
\item $p_1 = t$, $q_1 = \alpha$, $p_2 = y_{i,0} + \frac{1}{2} r^2$, $q_2 =
  \theta$, so that $\mu _i(t,\alpha,r,\theta) = (t,y_{i,0}+\frac{1}{2} r^2)$.
\item $\omega _i= dt \wedge d\alpha + r dr \wedge d\theta$.
\item $f _i(t,\alpha,r,\theta) = \frac{1}{2} r^2$.
\end{enumerate} 
Finally, to get the desired convexity, the radial vector field on
$\R^2$ radiating out from the origin lifts to a Liouville vector field
$V_i$ defined on $X_i \setminus \pi_i^{-1}(E_{i,0} \cup \ldots \cup
E_{i,n})$ which is, in $(p_1,q_1,p_2,q_2)$ coordinates, given by $V_i
= p_1 \partial_{p_1} + p_2 \partial_{p_2}$. (The fact that $V_i$ is
defined on $\pi_i^{-1}(E_{i,n_i+1})$ comes from the fact that the line
$L_{i,n_i+1}$ passes through $O$, so that $E_{i,n_i+1}$ is tangent to
the radial vector field.)  On the end
$\mu_i^{-1}((x_{i,0},x_{i,0}+\epsilon) \times
[y_{i,0},y_{i,0}+\epsilon)) \cong (x_{i,0},x_{i,0}+\epsilon) \times
S^1 \times D$, with coordinates $(t,\alpha,r,\theta)$ as above, $V_i$ is
given by $V = t \partial_t + (\frac{1}{2} r + \frac{y_{i,0}}{r})
\partial_r$.

Thus we see that $(X_i,\omega_i,C_i=\mu_i^{-1}(E_{i,0} \cup \ldots
\cup E_{i,n_i}), f_i = g_i \circ \mu_i, V_i)$ is the desired neighbourhood
$5$--tuple for the legs and that, furthermore, we have everything
written down explicitly in local coordinates on the end of $X_i$ which
is a neighbourhood of the disk $D_{i,0}$.

Now we construct the neighbourhood $5$--tuple for the central vertex
using the following lemma (which we have stated in a form which allows
also for positive genus, in case it should ever be useful):

\begin{lem} \label{L:centralvertex}
  Suppose that we are given a compact connected surface $\Sigma$ with $m>0$
  boundary components $\partial_1 \Sigma, \ldots, \partial_m \Sigma$, $k \geq
  0$ negative real numbers $c_1, \ldots, c_k$, and $m-k > 0$ positive real
  numbers $c_{k+1}, \ldots, c_m$, with $c_1 + \ldots + c_k + c_{k+1} + \ldots
  + c_m > 0$. Then there exists a symplectic form $\beta$ on $\Sigma$ and a
  Liouville vector field $W$ defined on all of $\Sigma$, pointing in along
  $\partial_i \Sigma$ for $1 \leq i \leq k$ and pointing out along $\partial_i
  \Sigma$ for $k+1 \leq i \leq m$, such that, for each $i = 1, \ldots, m$, a
  collar neighbourhood of $\partial_i \Sigma$ can be parameterized as
    $(c_i-\epsilon,c_i]\times S^1$  with coordinates $(t,\alpha)$
  with respect to which $\beta = dt \wedge d\alpha$ and $W = t \partial_t$.
\end{lem}

\begin{proof}
Note that the condition that $c_1 + \ldots + c_m > 0$ is a necessary
condition because $2\pi(c_1 + \ldots + c_m) = \int_{\partial \Sigma} t
d\alpha = \int_{\partial \Sigma}
\imath_W \beta = \int_\Sigma \beta > 0$.  There are probably numerous
ways to see that this lemma is true; in fact we need only extend the
$1$--forms $t d\alpha$ on each end to a $1$--form $\gamma$ on all of
$\Sigma$ such that $d\gamma > 0$. Our idea here is to apply
Weinstein's techniques~\cite{Weinstein} for attaching symplectic
handles to symplectic manifolds with convex boundary in the relatively
trivial case of dimension $2$; in this case the ``contact forms'' on
the $1$--dimensional boundary components are nowhere zero
$1$--forms. Once we get the $1$--forms on the boundary correct, then
flowing in along the Liouville vector fields produces the correct
parameterization of the ends.

First we claim that the lemma is true if $m-k=1$, in which
case the vector field should point out along only one boundary
component. Start with the disjoint union $\amalg_{i=1}^k ([c_i-\epsilon,c_i] \times S^1,
dt \wedge d\alpha)$ of $k$ symplectic cylinders with the Liouville
vector field $t \partial_t$
(or a very small disk if $k=0$ with the radial vector field) and
attach Weinstein's $1$--handles to the boundaries along which the
Liouville vector fields point out, to make a
symplectic surface diffeomorphic to $\Sigma$, with a Liouville vector
field with the correct behaviour along $\partial_1 \Sigma, \ldots,
\partial_k \Sigma$ and pointing out along $\partial_m \Sigma$. Now
simply attach a symplectization collar to $\partial_m \Sigma$ if
$\int_{\partial_m \Sigma} \imath_W \beta$ is too small, or remove a
collar (by flowing backwards along $W$) if it is too big, and thus
achieve the correct behaviour along $\partial_m \Sigma$.

Now if $m-k > 1$, we can first build a symplectic surface
$(\Sigma',\beta')$ with the same genus as $\Sigma$ and $k+1$ boundary
components, with a Liouville vector field $W$, with $k$ boundary
components along which $W$ points in and the behaviour is correct, and
one boundary component $\partial_+$ along which the vector field
points out, arranging that $\int_{\partial_+} \imath_{W'} \beta' <
2\pi(c_{k+1} + \ldots c_m)$. Then use the same ideas as in the preceding
paragraph to build a $(m-k+1)$--punctured sphere with one boundary
component $\partial_-$ along which the Liouville vector field points
in, so as to be able to glue to $\partial_+$, and with $m-k$
components along which the vector field points out with the correct
behaviour. Then glue the two pieces together.
\end{proof}

The neighbourhood $5$--tuple needed for the central vertex is now $(X_0 =
\Sigma \times D, \omega_0= \beta + \frac{1}{2\pi}r dr \wedge d\theta, C_0 =
\Sigma \times \{(0,0)\}, f_0 = \frac{1}{2} r^2, V_0 = W + (\frac{1}{2} r +
\frac{y_0}{r}) \partial_r)$. Because $y_0 = y_{i,0}$ for each $i=1, \ldots,
m$, if we implement the above lemma using $c_i = x_{i,0} + \delta$ for some
small $\delta > 0$, and for $\Sigma$ an $m$--punctured sphere, we will be able
to glue the $5$--tuples together. The area of the central sphere $S_0$ will
then be $2\pi(x_{1,1} + x_{2,1} + \ldots + x_{m,1})$ by Stokes' theorem.

It remains to compute $S_0 \cdot S_0$.
To do this note that, for some small $\delta > 0$, the submanifold $N$
of $X$ defined by:
\[ 
N = (\Sigma \times D_{\sqrt{2\delta}}) \cup \bigcup_{i=1}^m
\mu_i^{-1}( R_i \cap ((x_{i,0},\infty) \times
[y_{i,0},y_{i,0}+\delta))) 
\]
(where $D_{\sqrt{2\delta}}$ is an open disk of radius
  $\sqrt{2\delta}$) is a tubular neighbourhood of $S_0$. Then we can
  see $\partial N$ as diffeomorphic to $\Sigma \times S^1$ with the
  $i$'th component $(\partial_i \Sigma) \times S^1$ of $\partial
  (\Sigma \times S^1)$ filled in with a solid torus so that the
  $(-v_{i,1},u_{i,1})$ curves in $(\partial_i \Sigma) \times S^1 = S^1
  \times S^1$ bound disks. But since $v_{i,1} = 1$ this means that the
  $(1,-u_{i,1})$ curves are filled, so this $3$--manifold is the
  $S^1$--bundle over $S^2$ of Euler class $-u_{1,1} - u_{2,1} - \ldots
  - u_{m,1}$. Since the base $S^2$ in this case is precisely the
  central vertex sphere $S_0$, this tells us that $S_0 \cdot S_0 =
  -u_{1,1} - \ldots - u_{m,1}$.
\end{proof}

Now we need to translate Proposition~\ref{p:technicaltoric} into a
proof of Theorem~\ref{t:toric}. 

\begin{proof}[Proof of Theorem~\ref{t:toric}]

Given a star-shaped plumbing graph $\Gamma$, label the spheres $S_0,
S_{i,j}, 1 \leq i \leq m, 1 \leq j \leq n_i$, so that $S_0$
corresponds to the central $m$--valent vertex of $\Gamma$ and
$S_{i,1}, \ldots, S_{i,n_i}$ correspond in order to the vertices on
the $i$'th leg, with $S_0 \cdot S_{i,1} = 1$. Let the corresponding
self-intersections be $s_0, s_{i,j}$ and let the areas be $a_0,
a_{i,j}$. We need to produce the points $P_{i,j}, 1 \leq i \leq m, 0
\leq j \leq n_i+2$ satisfying the conditions of
Proposition~\ref{p:technicaltoric}, giving the desired
self-intersections and areas.

We first consider the vectors $\tau_{i,j} \in \Z^2$, $1 \leq i \leq m,
0 \leq j \leq n_i+1$. We are forced to have $\tau_{i,0} = (1,0)^T$ and
$\tau_{i,1} = (u_{i,1},1)^T$. We can make a choice for the values of
$u_{i,1}$ so long as $u_{1,1} + \ldots + u_{m,1} = -s_0$, so that $S_0
\cdot S_0 = - u_{i,1} - \ldots - u_{m,1} = s_0$. Now note that each
$\tau_{i,j+1}$ for $j=1, \ldots, n_i$ is completely determined by
$\tau_{i,j-1}$ and $\tau_{i,j}$ and the constraints that
$\det(\tau_{i,j-1},\tau_{i,j+1}) = -s_{i,j}$ and
$\det(\tau_{i,j},\tau_{i,j+1}) = +1$, so that our choices for
$\tau_{i,0}$ and $\tau_{i,1}$ determine $\tau_{i,2}, \ldots,
\tau_{i,n_i+1}$. In fact, $\tau_{i,j+1} = - \tau_{i,j-1} - s_{i,j}
\tau_{i,j}$.

For each $i = 1, \ldots, m$, let $\sigma_i = u_{i,n_i+1}/v_{i,n_i+1}$,
the reciprocal of the slope of the terminal vector
$\tau_{i,n_i+1}$. We claim that if $\Gamma$ is negative definite then,
for any choice of $u_{1,1}, u_{2,1}, \ldots, u_{m,1}$ such that
$u_{1,1} + u_{2,1} + \ldots + u_{m,1} = -s_0$, we will automatically
have $\sigma_1 + \ldots + \sigma_m > 0$. Given this claim, the
existence of suitable points $P_{i,j}$ can be seen as follows:

Having determined the vectors $\tau_{i,j}$, the points $P_{i,j}$ for
$i=1, \ldots,m, j = 2, \ldots, n_i+1$ are completed determined by the
$P_{i,1}$'s and the areas $a_{i,j}$. Furthermore, the fact that
$P_{i,n_i+1}$ must lie on the line $L_{i,n_i+1}$ passing through the
origin tangent to the vector $\tau_{i,n_i+1}$ constrains $P_{i,1}$ to
lie on a particular line parallel to $\tau_{i,n_i+1}$; i.e. there is
some constant $K_i$, determined by the $\tau_{i,j}$'s and the
$a_{i,j}$'s, such that $x_{i,1} = \sigma_i y_{i,1} + K_i$. Because the
points $P_{i,j}$ must all lie to the left of the line $L_{i,n_i+1}$,
which passes through the origin, we know that $K_i < 0$. The precise
location of $P_{i,n_i+2}$ on $L_{i,n_i+1}$ is not important, so we
ignore it. The precise location of $P_{i,0}$ is also not important,
provided it is close enough to $P_{i,1}$ and provided $x_{1,1} +
x_{2,1} + \ldots + x_{m,1} > 0$. Recall that we have the additional
constraint that $y_{i,1}=y_{i,0} = y_0$ for some $y_0$; thus we can
think of each $x_{i,1}$ as determined by a linear function of $y_0$:
$x_{i,1} =
\sigma_i y_0 + K_i$. We must now choose $y_0$ appropriately so that
$2\pi(x_{1,1}+x_{2,1}+ \ldots + x_{m,1}) = a_0$, the given area of
$S_0$. We have $x_{1,1} + x_{2,1} + \ldots + x_{m,1} =
(\sigma_1+\ldots+\sigma_m) y_0 + (K_1 + \ldots + K_m)$. But since
$\sigma_1 + \ldots + \sigma_m > 0$ and $K_1+\ldots+K_m < 0$, we can
realize any given positive value for $x_{1,1}+\ldots+x_{m,1}$ by
choosing an appropriate $y_0 > 0$.

Now to prove the claim, we first appeal to~\cite[Theorem~5.2]{NR}
stating that a star-shaped graph $G$ is negative definite if and only
if $s_0 + r_1 + \ldots + r_m < 0$, where $r_i$ is the negative
reciprocal of the continued fraction corresponding to the $i$'th leg,
i.e.:
\begin{equation}\label{e:frac}
r_i = - \frac{1}{s_{i,1} - \frac{1}{s_{i,2} - \ldots - \frac{1}{s_{i,n_i}}}} 
\end{equation}
In the proof of \cite[Theorem~9.20]{Sym4from2} Symington shows that if
$\tau_0, \ldots, \tau_{n+1}$, $\tau_i = (u_i,v_i)^T$, is a list of vectors in
$\Z^2$ with $\tau_0 = (0,-1)^T$, $\tau_1 = (1,0)^T$, $\det(\tau_i,\tau_{i+1})
= 1$ and $\det(\tau_{i-1}, \tau_{i+1}) = b_i$, then the reciprocal slope
$\sigma$ of the last vector is given by:
\[ 
\sigma = \frac{u_{n+1}}{v_{n+1}} = b_1 - \frac{1}{b_2 - \ldots
  -\frac{1}{b_n}} 
\] 
We can change such a list to one where $\tau_0 = (1,0)^T$ and $\tau_1
= (u,1)^T$ for some given $u \in \Z$ by applying the linear
transformation
\[ \tau_i \mapsto \left( \begin{array}{c c} u & -1 \\ 1 & 0
\end{array} \right) \tau_i \]
This transforms $\sigma$ to $u - 1/\sigma$. In other words, after
applying this transformation, $\sigma$ is now given by:
\[ \sigma = u - \frac{1}{b_1 - \frac{1}{b_2 - \ldots -
    \frac{1}{b_n}}} \]
Applying this to each leg of our star-shaped graph, with the initial
condition $\tau_{i,1} = (u_{i,1},1)^T$ and
$\det(\tau_{i,j-1},\tau_{i,j+1}) = -s_{i,j}$ we get that 
\[ \sigma_i = u_{i,1} + \frac{1}{s_{i,1} - \frac{1}{s_{i,2} - \ldots -
    \frac{1}{s_{i,n_i}}}} = u_{i,1} - r_i \]
Thus $\sigma_1 + \ldots + \sigma_m = -(s_0 + r_1 + \ldots + r_m)$ and
the claim is proved.
\end{proof}

\begin{rem} {\rm Since Lemma~\ref{L:centralvertex} allows for any genus,
    Theorem~\ref{t:nbhood} is actually true even when the surface
    corresponding to the central vertex has positive genus.  It seems that a
    similar picture provides $\omega$--convex neighbourhood for any collection
    of symplectic spheres intersecting each other according to a
    \emph{negative definite} plumbing tree, i.e. the assumption of
    Theorem~\ref{t:toric} on $\Gamma$ being star-shaped can be removed. We
    hope to return to this question in a future project.}
\end{rem}

\section{Tight contact structures on certain Seifert fibered $3$--manifolds}

In order to find the right gluing map for performing the rational blow-down
process, we will invoke a classification result of tight contact structures on
certain small Seifert fibered $3$--manifolds. Let us start with some
generalities.

\begin{figure}[ht]
\begin{center}
\psfrag{0}{$s_0$}
\psfrag{1}{$-\frac{1}{r_1}$}
\psfrag{2}{$-\frac{1}{r_2}$}
\psfrag{3}{$-\frac{1}{r_3}$}
\includegraphics[width=10cm]{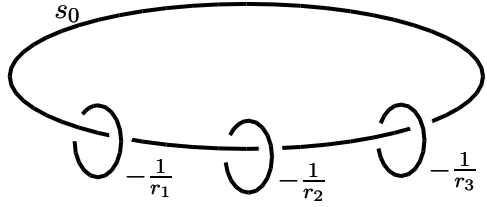}
\end{center}
\caption{Surgery diagram for the Seifert fibered $3$--manifold
$M(s_0;r_1,r_2,r_3)$}
\label{f:seif}
\end{figure}

\begin{figure}[ht]
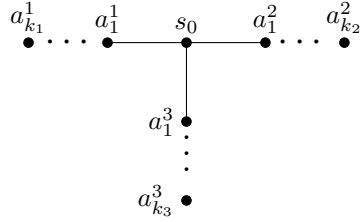

\begin{center}
\setlength{\unitlength}{1mm}
\unitlength=0.7cm
\begin{graph}(10,4)(0,-3)
\graphnodesize{0.2}

 \roundnode{m2}(2,0)
 \roundnode{m3}(3.5,0)
 \roundnode{m4}(5,0)
 \roundnode{m5}(6.5,0)
 \roundnode{m6}(8,0)
 \roundnode{m8}(5,-1.5)  
 \roundnode{m9}(5,-3)

\edge{m3}{m4}
\edge{m4}{m5}
\edge{m4}{m8}

  \autonodetext{m2}[n]{{\small $a^1_{k_1}$}}
  \autonodetext{m3}[n]{{\small $a_1^1$}}
  \autonodetext{m4}[n]{{\small $s_0$}}
  \autonodetext{m5}[n]{{\small $a^2_1$}}
  \autonodetext{m6}[n]{{\small $a^2_{k_2}$}}
  \autonodetext{m8}[w]{{\small $a^3_1$}}
  \autonodetext{m9}[w]{{\small $a^3_{k_3}$}}
  \autonodetext{m3}[w]{{\Large $\cdots$}}
  \autonodetext{m5}[e]{{\Large $\cdots$}}

\freetext(5,-1.8){\Large $.$}
\freetext(5,-2.1){\Large $.$}
\freetext(5,-2.4){\Large $.$}

\end{graph}
\end{center}
\caption{Plumbing diagram of a $4$--manifold with boundary $M(s_0;r_1,r_2,r_3)$}
\label{f:plumbi}
\end{figure}

As is costumary, we say that a $3$--manifold $Y$ is a \emph{small
Seifert fibered space} if it can be given by the surgery diagram of
Figure~\ref{f:seif}. Here we  assume that $s_0\in \bfz $ and $r_i
\in (0,1)\cap \bfq$ (then $(s_0; r_1, r_2, r_3)$ are the
\emph{normalized Seifert invariants} of $M=M(s_0; r_1, r_2, r_3$)).
By applying the inverse slam dunk operation
(cf. \cite[Figure~5.30]{GS}), the diagram of Figure~\ref{f:seif} can
be easily transformed into a star-shaped plumbing as it is given by 
Figure~\ref{f:plumbi}.
The plumbing coefficients $[a^i_1, \ldots , a^i _{k_i}]$ on the $i^{th}$ leg
(satisfying $a^i_j \in \bfz $ and $a^i_j \leq -2$) are specified by the
continued fraction coefficients of the rational number $-\frac{1}{r_i}<-1$
(cf. Equation~\eqref{e:frac}). Notice that all plumbing graphs in the
Introduction (Figures~\ref{f:wahltype}--\ref{f:cs}) give rise to
$4$--manifolds with small Seifert fibered $3$--manifold boundary.

In a remarkable series of papers \cite{OSzF1, OSzF2} Ozsv\'ath and Szab\'o
introduced an invariant for spin$^c$ $3$--manifolds: the (mod 2) 
\emph{Ozsv\'ath--Szab\'o homology group} $\hf (Y, \t )$ of a closed 
spin$^c$ $3$--manifold $(Y, \t )$ is a finite dimensional vector space over the
field $\bfz _2$ of two elements. For a rational homology sphere the
dimension of this vector space is odd for every spin$^c$ structure. 
We say that a rational homology sphere $Y$ is an \emph{$L$--space}
if $\hf (Y , \t )=\bfz _2$ for every spin$^c$ structure $\t\in Spin^c(Y)$, 
that is, $Y$ has the simplest possible Ozsv\'ath--Szab\'o homologies.

Recall that a contact structure $\xi$ on any $3$--manifold naturally
induces a spin$^c$ structure, which we will denote by $\t _{\xi }$. In
addition, $(Y, \xi )$ also gives rise to an element $c(Y, \xi )\in \hf
(-Y, \t _{\xi })$, the \emph{contact invariant} of $(Y, \xi )$.
(For the definition and basic properties of $c(Y, \xi )$ see 
\cite{OSzcont, LSI}.)

As the next result (a compilation of theorems of Wu, Ghiggini and 
Plamenevskaya) shows, tight contact structures on small Seifert fibered
$3$--manifolds with simple Ozsv\'ath--Szab\'o homologies (and $s_0\leq -2$)
admit a simple classification scheme. More precisely

\begin{thm}\label{t:legf3}
Suppose that the small Seifert fibered $3$--manifold $M=M(s_0; r_1, r_2,
r_3)$ satisfies $s_0\leq -2$ and $M$ is an $L$--space.  Then two tight
contact structures $\xi _1, \xi _2$ on $M$ are isotopic if and only if
$\t _{\xi _1}=\t _{\xi _2}$.
\end{thm}

\begin{proof}
One direction of the statement is obvious: isotopic contact structures
are homotopic as $2$--plane fields, hence induce the same spin$^c$
structure.  The converse, however, is more subtle, since it states
that in these $3$--manifolds the spin$^c$ structure (which is a
homotopic invariant of the contact structure) determines the tight
contact structure up to isotopy.

Let $\Gamma$ denote the star-shaped plumbing graph we get from the
surgery diagram of $M=M(s_0; r_1, r_2, r_3)$ by inverse slam
dunks. The corresponding $4$--manifold $Z_{\Gamma}$ is then given by a
sequence of $2$--handle attachments along the unknots corresponding to
the vertices of $\Gamma$. Working with a projection now it is easy to
see that each unknot can be isotoped until it becomes the Legendrian
unknot, i.e., a Legendrian knot isotopic to the unknot with
Thurston--Bennequin invariant $-1$.  The assumption $s_0\leq -2$
ensures that by adding sufficiently many zig-zags to these Legendrian
unknots we end up with a Legendrian link on which contact
$(-1)$--surgery results $Z_{\Gamma}$, together with a Stein structure
on it \cite{Eli2, Gompf}. Notice, however, that there is a certain freedom in
adding the zig-zags to the unknots: moving a zig-zag from left to
right (or from right to left) will change the rotation number of the
particular knot, which in turn will change the first Chern class of
the resulting Stein structure on $Z_{\Gamma}$.  According to a result
of Lisca and Mati\'c \cite{LM} Stein structures on a fixed $4$--manifold
with different first Chern classes induce nonisotopic tight contact
structures on the boundary of the $4$--manifold. The number of tight
contact structures one can distinguish in this way can be easily
computed from the framing coefficients of the graph $\Gamma$.
Furthermore, it was shown by Plamenevskaya~\cite{Olga} that if two
Stein structures on a $4$--manifold induce two contact structures $\xi
_1, \xi _2$, and for the two Stein structures we have $c_1(J_1)\neq
c_1(J_2)$ then for the (mod 2 reduced) contact Ozsv\'ath--Szab\'o
invariants we have $c(M, \xi _1)\neq c(M, \xi _2)$.  In summary, the
Legendrian surgery construction described above allowed us to
construct a finite set of Stein fillable (hence tight) contact
structures on $M$ (corresponding to the different choices of left and
right zig-zags), and all these structures have different (and nonzero)
contact Ozsv\'ath--Szab\'o invariants.

Now we are ready to prove the theorem. Using convex surface techniques, it was
shown by Wu \cite{Wu} that for a small Seifert $3$--manifold with $s_0\leq -3$
the above set of Stein fillable contact structures, in fact, contains
\emph{all} tight contact structures on the $3$--manifold at hand. (For these
$3$--manifolds the $L$--space condition is automatically satisfied, cf.
Lemma~\ref{l:sat}.) For the remaining $s_0=-2$ case we argue as follows:
According to~\cite{LSIII} the assumption for $M$ being an $L$--space implies
that $-M$ admits no transverse contact structure. In~\cite{ghigg} this
property is translated to a numerical condition for $r_1,r_2,r_3$, which in
turn, implies that the number of tight contact structures on $M$ is bounded
above by the number of Stein fillable contact structures constructed by
Legendrian surgery along the Legendrian links described above.  In short,
under the assumption of the theorem, the Legendrian surgery construction
described at the beginning of the proof produces all tight contact structures
on $M$.

Appealing to the $L$--space property again, now we can easily conclude
the proof: if $\xi _1$, $\xi _2$ are nonisotopic tight contact
structures on $M$ then according to the above said we have $c(M, \xi
_1)\neq c(M, \xi _2)$. Assuming that for the induced spin$^c$
structures $\t _{\xi _1}=\t _{\xi _2}=\t $ holds we easily reach a
contradiction: in this case the contact invariants $c(M, \xi _1), c(M,
\xi _2)\in \hf (-M, \t )=\bfz _2$ are distinct, nonzero elements,
which is clearly impossible.
\end{proof}

\begin{rem} {\rm In fact, a very similar proof (together with the
    classification of tight contact structures on small Seifert fibered
    $3$--manifolds with $s_0\geq 0$~\cite{GLS}) applies for $s_0\geq 0$ as
    well.  (It is worth mentioning that if $s_0(M)\geq 0$ then $s_0(-M)\leq
    -3$, therefore any $M$ with $s_0(M)\geq 0$ is still an $L$--space.)  The
    only difference is in the construction of the appropriate $4$--manifold
    (carrying the Stein structures): for $s_0\geq 0$ we need to introduce an
    appropriate Stein $1$--handle as well, cf. the argument of~\cite{GLS}.  We
    will not use this fact in the present paper. A similar statement (that is,
    that tight contact structures with isomorphic spin$^c$ structures are
    isotopic) is expected in the case when $s_0=-1$ and the $3$--manifold is
    an $L$--space; for related discussion and partial results see~\cite{GLS2}.
    Most probably the same statement also holds for \emph{strongly fillable}
    contact structures on boundaries of negative definite plumbings which are
    $L$--spaces.  Considering all tight contact structures, such a statement
    cannot be true, since many of these $3$--manifolds are toroidal, hence
    contain infinite families of homotopic, nonisomorphic tight contact
    structures.  (These infinite families are constructed by inserting Giroux
    torsions, resulting in contact structures which are not strongly
    fillable~\cite{gay}.)  For strongly fillable structures, however, one can
    expect that the strong filling can be chosen to be diffeomorphic to the
    plumbing $4$--manifold, and then the adaptation of the proof above would
    provide the same result.}
\end{rem}

\section{The proof of Theorem~\ref{t:main}}

We start this section by proving that the result of
Theorem~\ref{t:legf3} applies for the plumbing manifolds given by
graphs in $\frg$.

\begin{lem}\label{l:sat}
Suppose that $\Gamma \in \frg$. Then the star-shaped graph $\Gamma $ is
negative definite and the $3$--manifold $\partial M_{\Gamma}$ is an
$L$--space.
\end{lem}
\begin{proof}
The plumbing graph $\Gamma \in \frg$ is negative definite since it
embeds in a negative definite lattice \cite{SSW}. (Alternatively, a
simple computation and the application of \cite[Theorem~5.2]{NR} shows
the same.) There are many ways to verify that $\partial M_{\Gamma}$ is
an $L$--space. The direct application of the algorithm
of~\cite{OSzplum} (which applies for negative definite plumbing graphs
with at most one ``bad'' vertex) implies the result at once, although
this proof involves some computations. (In fact, for graphs in
$\frw\cup \frn$ the result is stated in~\cite{OSzplum}, since these
graphs involve no ``bad'' vertices.)  Alternatively, we can argue as
follows: as~\cite[Examples~8.3 and 8.4]{SSW} show, the normal surface
singularity defined by $\Gamma\in \frg$ admits a rational homology disk
smoothing and this property implies that the surface singularity is
rational, cf.~\cite[Proposition~2.3]{SSW}. (In yet another way,
Laufer's algorithm can be easily applied to verify rationality of the
surface singularities defined by graphs in $\frg$.)  The computation
presented in~\cite{nemethi} shows that the link of a rational
singularity is an $L$--space, concluding our argument.
\end{proof}

After all these preparations, now we are ready to turn to the proof of 
Theorem~\ref{t:main}:

\begin{proof}[Proof of Theorem~\ref{t:main}]
  Suppose that $S_i$ are symplectic spheres in the sym\-plec\-tic
  $4$--manifold $(X, \omega )$, intersecting each other
  $\omega$--orthogonally, and according to the plumbing tree $\Gamma
  \in \frg$. By Theorem~\ref{t:nbhood} there is a neighbourhood
  $S_{\Gamma }$ of $\cup S_i$ such that its complement is a strong
  concave filling of its contact boundary. The contact structure on
  the boundary $\partial S_{\Gamma}$ induces the spin$^c$ structure
  which extends to $S_{\Gamma}$ as $\s$ with the property that $c_1(\s
  )$ satisfies the adjunction equality for all $S_i$ (since $S_i$ are
  all smooth symplectic submanifolds). 

Now consider the normal surface singularity given by $\Gamma$. As
  such, it defines a contact structure on its link, which induces the
  spin$^c$ structure $\t$.  Resolving the singularity doesn't change
  this contact structure (and hence the spin$^c$ structure $\t$), but
  provides a $4$--manifold diffeomorphic to $S_{\Gamma}$, with a
  spin$^c$ structure $\s _{can}$.  Since the vertices of $\Gamma$ in
  the resolution correspond to complex curves, the first Chern class
  of $\s _{can}$ also satisfies the adjunction equality. This shows
  that $c_1(\s )=c_1(\s _{can})$, and since $S_{\Gamma }$ is a simply
  connected $4$--manifold, we conclude that $\s =\s _{can}$. This fact,
  however, identifies the spin$^c$ structures of the contact
  structures on $\partial S_{\Gamma }$, which in the light of
  Theorem~\ref{t:legf3} implies that the contact structure on
  $\partial S_{\Gamma }$ induced by the Liouville vector field is
  isotopic to the one induced on the link of the singularity.  
This last observation verifies the existence of the desired gluing
  contactomorphism, hence the symplectic gluing can be
  performed as in~\cite{etny}, concluding the proof.
\end{proof}

\end{document}